\documentclass{amsart}
\usepackage{graphicx}
\usepackage{hyperref}
\usepackage{amssymb}
\usepackage{tikz-cd} 
\usepackage{amsmath,amscd,amsthm}
\usepackage[mathscr]{euscript}
\usepackage[all]{xy}
\usepackage{lmodern}
\usepackage[T1]{fontenc}
\usepackage{tikz}
\usetikzlibrary{arrows,automata}
\usetikzlibrary{matrix,arrows,decorations.pathmorphing}
 \usepackage{color}
\DeclareMathOperator{\Surj}{Surj}

\DeclareMathOperator{\Orb}{Orb}
\newcommand{\bbZ}{\mathbb Z}
\newcommand{\bbN}{\mathbb N} 
\newcommand{\bbF}{\mathbb F} 
\newcommand{\cF}{\mathcal F}
\newcommand{\cX}{\mathcal X}
 
\newtheorem{theorem}{Theorem}[section]
\newtheorem{lemma}[theorem]{Lemma}
\newtheorem{proposition}[theorem]{Proposition}
\newtheorem{corollary}[theorem]{Corollary}

\theoremstyle{definition}

\newtheorem{example}[theorem]{Example}
\newtheorem{remark}[theorem]{Remark}
\begin{document}

\title{Identifying domatic partitions via graph dynamical systems}
\author{Mehmet Akif Erdal}
\address{Faculty of Arts and Sciences, Department of Mathematics,
Yeditepe University, 34755, İstanbul, Turkiye}
\email{mehmet.erdal@yeditepe.edu.tr}
 
\begin{abstract}
	For a graph $G=(V,E),$ a set of  vertices $D\subseteq V $ is called a dominating set if every vertex  in $V\backslash D$ is adjacent to a vertex in $D.$ A domatic-$2$-partition of $G$ is a partition of its vertices into two disjoint dominating sets. In this paper, for a finite simple connected graph $G,$ we construct a graph dynamical system $F$ and show that the set of dominating sets of $G$ are in one-to-one correspondence with the image of the action map of $F$. Moreover, we obtain the set of all domatic-$2$-partitions of $G$ from the set of all periodic orbits of $F.$ Finally, we extended actions of two dynamical systems  to an action of a free semigroup on two letters, and determine independent dominating sets and idomatic partitions using its maximal invariant subset with a reversible action.
\end{abstract}	
\keywords{graph dynamical system, dominating set, domatic partition}
\maketitle

	\section{Introduction}
	For a given graph $G=(V,E),$ a graph dynamical system on $G$ consists of a set $S$ of states and an update function $F:S^V\to S^V$ on the set of functions from $V$ to $S,$ where $F$ is defined via a local update functions $f_v:S^{N[v]}\to S$ for all  $v\in V.$  Here $N[v]$ denotes the closed neighborhood of $v$. This notion was first introduced in \cite{barrett2000elements}. Later many applications on modeling and simulation of complex systems, such as biological and social networks, have been studied, see e.g. \cite{eubank2004modelling,barrett2011modeling,saadatpour2010attractor}. If  the state set is $\{0,1\}$ and the update sequence is synchronous then such a graph dynamical system is often called a parallel Boolean dynamical system. In this paper, we use such systems but for simplicity we just refer them as `graph dynamical system'.

For any set $X$, we denote by $\bbF_2^{X}$ the set of functions from $X$ to $\bbF_2$.
Let $G=(V,E)$  be a finite simple connected graph. We define a graph dynamical system   $F:\bbF_2^{V}\to \bbF_2^{V}$ by $$F({s})={s} P_{s} +P_{s} + {s} +1$$ where $$P_{s}(v)=\prod_{x\in N(v)}(1+s(x)) $$ for $v\in V.$ This graph dynamical system has many interesting properties; e.g., see Lemmas \ref{lem:preimageoffsigma}, \ref{lem:Fs=F(1+Fs)} or \ref{lem:Fs=1+s}. We show that the set of dominating sets of the graph $G$ are in one-to-one correspondence with configurations in the image of $F,$ see Corollary \ref{cor:inverseimageofzeroisdomset}. The correspondence is given by sending the dominating set to the function that sends every vertex of the dominating set to $0$ and other vertices to $1$.  Besides, this system has only one fixed point; namely, the constant function $c_0$ at $0,$ and has non-zero periodic points. We also show that all of these periodic points  have periods $2$ and pairs of such points are in one-to-one correspondence with  domatic-$2$-partitions of $G,$ see Theorem \ref{thm:main}.

The paper is organized as follows. In Section \ref{sec:prelim} we give some background on dominating sets, domatic partitions and graph dynamical systems. In Section  \ref{sec:main} we state and prove our main results; particularly, Corollary \ref{cor:inverseimageofzeroisdomset} and Theorem \ref{thm:main}.
Finally in Section \ref{ssect:moredominating} we extended our system to two different semigroup dynamical systems  and present their applications on independent dominating sets and idomatic partitions (that is, domatic partitions by independent dominating sets) of graphs, see Corollaries \ref{cor:maxindep} and \ref{cor:idomatic}.

\section{A few preliminaries}\label{sec:prelim}
Throughout this paper  $G=(V,E)$ denotes an undirected finite simple connected graph with vertex set $V$ and edge set $E.$
\subsection{Dominating sets and domatic partitions}
Dominating sets of graphs and related problems are well studied in the literature. For general references about classical and recent advances we refer to \cite{haynes1998fund} and \cite{haynes2020topic}. A dominating set for $G$ is a subset $D$ of $V$ for which every vertex that is not in $D$ is adjacent to a vertex in $D.$  

A domatic partition for $G$ is a partition of its vertex set $V$ into dominating sets. If the number of dominating sets in the partition is $n,$ then it is called domatic-$n$-partition of $G.$ Thus, by definition if $V$ is written as a disjoint union as $V=D\sqcup B$ where $D$ and $B$ are dominating sets, then the partition $\{D,B\}$ is a domatic $2$-partition of $G.$

\subsection{Graph dynamical systems}
For a vertex $v\in V $ the open  neighborhood of $v$ is given by $$N(v)=\{x\in V \ \mid  \ {vx}\in E \}$$ and the closed  neighborhood of $v$ is $N[v]=N(v) \cup \{v\}.$ A graph dynamical system consists of a graph $G=(V,E),$ a finite state set $S$ and an update function $F:S^V\to S^V$ where $F$ is defined via a local update function $f_v:S^{N[v]}\to S.$ For more details on graph dynamical systems we refer to \cite{barrett2000elements}. If $F$ satisfies  $F(s)(v)=f_v(s\mid _{N[v]})$ for each $v$ in $V,$ then such a system is referred as generalized cellular automata (see also \cite{sutner1989cellular}).  In this paper, we are interested in a special case where  $S=\bbF_2$  the field with $2$ elements $\{0,1\}.$

\section{Main Results}\label{sec:main}
For any set $X$ denote by $\bbF_2^X$ the set of functions from $X$ to $\bbF_2.$ The addition and multiplication of configurations in $\bbF_2^X$ is defined point-wise. For each ${s}:V\to \bbF_2$ let $P_{s}:V\to \bbF_2$ be the function given by $$P_{s}(v)=\prod_{x\in N(v)}(1+s(x)).$$ We define a graph dynamical system with update function $F:\bbF_2^{V}\to \bbF_2^{V}$ given by $$F({s})=1+ {s} +P_{s}+ {s} P_{s} =(1+s)(1+P_s).$$ In fact, this function  induces a generalized cellular automota with  local functions $f_v:\bbF_2^{N[v]}\to \bbF_2$ given by $$	f_v(s)=1+s(v)+\prod_{x\in N(v)}(1+s(x)) +(s(v)\prod_{x\in N(v)}(1+s(x))).$$ 
\begin{lemma}\label{lem:preimageoffsigma}
	For any $s\in \bbF_2^{V},$  $$F({s})^{-1}(1)=\{v\in {s}^{-1}(0)\ \mid  \ \exists\ u\in{s}^{-1}(1),  \ {uv} \in E\}.$$
\end{lemma}	
\begin{proof}
	Assume that $v\in F({s})^{-1}(1),$ so that $(1 + s)(v)=1=(1 + P_s)(v)$. Then  $${s}(v)=0=\prod_{x\in N(v)}(1+{s}(x)).$$ But this means there is $u\in  N(v)$ with ${s}(u)=1.$ Therefore, $u\in{s}^{-1}(1)$ and ${uv} \in E.$ Conversely, if ${s}(v)=0$ and there exists $u\in{s}^{-1}(1)$ such that ${uv} \in E,$ then $u\in N(v)$ and $1+{s}(u)=0,$ which implies $P_s(v)=0.$ 
\end{proof}
\begin{lemma}\label{lem:Fs=F(1+Fs)}
	For every $s\in \bbF_2^{V},$  $F(s)=F(1+F(s)).$ 
\end{lemma}	
\begin{proof}
	We have $$F(1+F(s))(v)=F(s)(v)(1+\prod_{x\in N(v)}F(s)(x)).$$ Then the equality holds whenever $F(s)(v)=0.$  If $s(v)=1,$ then $$F(s)(v)=0=F(1+F(s))(v).$$ If $s(v)=0,$ then $F(s)(v)=1$ if and only if there exists $u\in N(v)$ with $s(u)=1.$ But this implies that $$\prod_{x\in N(v)}F(s)(x)=0$$ as $F(s)(u)=0$ and $u\in N(v).$ Therefore $F(s)(v)=F(1+F(s))(v).$
\end{proof}

\begin{lemma}\label{lem:Fs=1+s}
	For any $s\in \bbF_2^{V},$ $F(s)=1+s$ if and only if  ${s}^{-1}(1)$ is a dominating set.
\end{lemma}	
\begin{proof}
	Assume that $F(s)=1+s.$ If $v\in s^{-1}(0),$ so that $F(s)(v)=1,$ then by lemma \ref{lem:preimageoffsigma} there exists $u\in V$ with $s(u)=1$ and ${uv} \in E.$ Thus, ${s}^{-1}(1)$ is a dominating set. For the converse, assume that ${s}^{-1}(1)$ is a dominating set. Then for every $v\in s^{-1}(0),$ there exists $u\in  {s}^{-1}(1)$  with ${uv}\in E$; i.e., there exists $u\in N(v)$ with $s(u)=1.$ This means $P_s(v)=0$ and $F(s)(v)=1.$ Note that if $s(v)=1,$ then $F(s)(v)=0.$ Therefore, for every $v\in V,$ $F(s)(v)=1+s(v).$
\end{proof}

\begin{corollary}\label{cor:inverseimageofzeroisdomset} 
	For any $s\in \bbF_2^{V},$ $F(s)^{-1}(0)$ is a dominating set. Conversely, for every dominating set $D\subseteq V,$ there exists a configuration $s\in \bbF_2^{V}$ such that $F(s)^{-1}(0)=D.$ In particular, there is a bijection between set of configurations in the image of $F$ the set of all  dominating sets of $G.$
\end{corollary}	
\begin{proof}
By lemma \ref{lem:Fs=F(1+Fs)} for every $s\in \bbF_2^V$, we have $F(1+F(s))=F(s)=1+(1+F(s))$. Thus, by lemma \ref{lem:Fs=1+s} we get that $(1+F(s))^{-1}(1)=F(s)^{-1}(0)$ is a dominating set. Conversely, if $D$ is a dominating set, define $s:V\to \bbF_2$ by $$s(v)=\begin{cases} 1 & \text{ if } v\in D \\
	0 & \text{ if } v\notin D 
\end{cases}$$
Then $F(s)=1+s$ and $s^{-1}(1)=F(s)^{-1}(0)=D$.
\end{proof}

Let ${\Surj}{(\bbF_2^{V})}= {\bbF_2^{V}}\backslash \{c_0,c_1\}$ where $c_i$ is the constant function at $i\in \{0,1\}$; i.e., the subset of $\bbF_2^{V}$ of all surjections from $V$ to $\bbF_2.$ A subset  $A$ of an $\mathbb N$-set is called  forward invariant  if for every $n\in N,$ $n\cdot A\subseteq A.$ We have $F(c_0)=c_0$ and $F(c_1)=c_0,$ so that $\{c_0,c_1\}$ is forward invariant under the action of $F$; i.e.,  $F(\{c_0,c_1\})\subseteq  \{c_0,c_1\}$.  We also have the following proposition.
\begin{proposition}\label{prop:surjisforwardinvariant} $F({\Surj}{(\bbF_2^{V})})\subseteq 	 {\Surj}{(\bbF_2^{V})}$.
\end{proposition}	
\begin{proof} Let ${s}\in {\Surj}{(\bbF_2^{V})}$ such that $F({s})\in \{c_0,c_1\}.$ Since   $F(s)^{-1}(0)$ is a dominating set for every $s\in \bbF_2^V$,   we have $F({s})\neq c_1.$ Assume $F({s})= c_0.$ Then, by lemma \ref{lem:preimageoffsigma}, $$\emptyset=F({s})^{-1}(1)=\{v\in {s}^{-1}(0)\ \mid  \ \exists\ u\in{s}^{-1}(1),  \ {uv} \in E\}.$$ But as $G$ is connected, there exist two adjacent vertices $v,w\in V$ such that ${s}(v)=1$ and ${s}(w)=0,$ so we get a contradiction. Thus, $$F({\Surj}{(\bbF_2^{V})})\subseteq 	 {\Surj}{(\bbF_2^{V})}.$$ \end{proof}
\begin{proposition}\label{prop:surjhasnofixedpoints}
	The only fixed point of $F$ is $c_0.$ In particular, ${\Surj}{(\bbF_2^{V})}$ does not have any fixed points under $F.$
\end{proposition}	
\begin{proof}
	First observe that $F(c_0)=c_0$. Assume that $s(v)=1$ for some $v\in V.$ Then  $$F(s)(v)={s}(v)P_{s}(v)  +P_{s}(v)+s(v) +1=0.$$ Thus $F(s)\neq s$; i.e., the only fixed point is $c_0$.
\end{proof}
Now we have the following theorem.
\begin{theorem}\label{thm:main} Let $G=(V,E)$ be a finite simple connected graph and let $F:\bbF_2^{V}\to \bbF_2^{V}$ be the system given by $$F({s})={s} P_{s} +P_{s} + {s} +1$$ where $P_{s}(v)=\prod_{x\in N(v)}(1+s(x)) $ for $v\in V.$ Then the set of all periodic points in ${\Surj}{(\bbF_2^{V})}$ is nonempty and every periodic point  has period $2.$ Moreover, the orbit space of these periodic points under the involution induced by $F$ is in one-to-one correspondence with the set of all domatic-$2$-partitions of $G.$ 
	
\end{theorem}	
\begin{proof}
	By Propositions \ref{prop:surjisforwardinvariant} and \ref{prop:surjhasnofixedpoints}, ${\Surj}{(\bbF_2^{V})}$ contains periodic points that are not fixed. Suppose that $s\in {\Surj}{(\bbF_2^{V})}$ is a periodic point of $F.$ Then there exist a positive integer $n$ such that $F^n(s)=s.$ This means $s$ is in the image of $F$; and thus, $s^{-1}(0)$ is a dominating set by Corollary \ref{cor:inverseimageofzeroisdomset}. Suppose that the period of $s$ is greater than $2;$ i.e.,  $F^2(s)\neq s.$ Then there exists $v\in V$ such that either $s(v)=1$ and $F^2(s)(v)=0$ or $s(v)=0$ and $F^2(s)(v)=1.$ 
	
	Assume first $s(v)=1$ and $F^2(s)(v)=0.$ Since $s(v)=1,$ we have $F(s)(v)=0.$ Then $$0=F^2(s)(v)=(1+F(s)(v))(1+ \prod_{x\in N(v)}(1+F(s)(x)))=1+ \prod_{x\in N(v)}(1+F(s)(x)).$$ Hence, we obtain $$\prod_{x\in N(v)}(1+F(s)(x))=1,$$ which implies that $F(s)(x)=0$ for every $x\in N(v).$ As ${ux}\in E,$ $v\in N(x),$ and as $s(v)=1$ we have $$\prod_{u\in N(x)}(1+s(u))=0.$$ Therefore, for every $x\in N(v),$ we must have  $s(x)=1.$ Since $s(v)=1$ and $s$ maps all of the neighboring vertices to $1,$ the vertex $v$ is not adjacent to  any vertex in $s^{-1}(0).$ But  this is a contradiction as  $s^{-1}(0)$ is a dominating set.
	
	Now assume that $s(v)=0$ and $F^2(s)(v)=1.$ Note that  $F(s)(v)=0$ (otherwise we should have $F(F(s))(v)=F^2(s)(v)=0$). Let $k$ be the greatest integer in $\{2, 3, \dots,  n-1\}$ such that $F^k(s)(v)=1.$ Denote the configuration $F^k(s)$ by $t.$ Then $t$ is a periodic point satisfying $t(v)=1$ and $F^2(t)(v)=0.$ But this reduces to the first case, which leads to a contradiction. Therefore, $n=2$; i.e., every periodic point has period $2.$ 
	
	Now let $s$ be a periodic point in ${\Surj}{(\bbF_2^{V})},$ so that $s$ has period $2.$ Assume that $F(s)\neq 1+s.$ Then there exists $v\in V$ such that $$F(s)(v)=0=s(v);$$ that is,  $$ \prod_{x\in N(v)}(1+(s)(x))=1.$$ Thus, we obtain $s(x)=0$ for all $x\in N(v).$ Assume that $F(s)(u)=1$ for some $u\in N(v).$ Then, $$\prod_{x\in N(v)}(1+F(s)(x))=0$$ and $$F^2(v)=(1+F(s)(v))(1+ \prod_{x\in N(v)}(1+F(s)(x)))=1,$$ which is a contradiction. Therefore, for all $x\in N(v)$ we have $F(s)(x)=0.$ So we have for all $x\in N(v),$ $$F(s)(x)=0=s(x).$$ By  repeating the same process for all elements in $N(x),$ we obtain that $$F(s)(u)=0=s(u)$$ for all $u\in N(x).$ As $G$ is a finite simple connected graph, applying the same process inductively, we obtain that $s(x)=0$ for every $x\in V$; i.e., $s=c_0.$ But this is a contradiction, since $c_0\notin {\Surj}{(\bbF_2^{V})}$. Therefore $F(s)=1+s$ and by lemma \ref{lem:Fs=1+s}, $s^{-1}(1)$ is a dominating set. Since we already have $s^{-1}(0)=V\backslash s^{-1}(1)$ is a dominating set, the pair  $\{s^{-1}(1), s^{-1}(0)\}$ forms a domatic-$2$-partition of $G.$

	Conversely, suppose that $\{D, B\}$ is a domatic-$2$-partition of $G.$ Define $s:V\to \bbF_2$ by $s(v)=1$ if $v\in D$ and $s(v)=0$ if $v\in B$; i.e., both ${s}^{-1}(0)$ and ${s}^{-1}(1)$ are dominating sets. By lemma \ref{lem:Fs=1+s} and Corollary \ref{cor:inverseimageofzeroisdomset}, $F(s)=1+s$ and $F(1+s)=s$; i.e., $s$ is a periodic point with period $2.$ The one-to-one correspondence between orbit space and set of domatic partitions is given by $\{s,1+s\}\mapsto \{s^{-1}(1), s^{-1}(0)\}.$ 
	
\end{proof}
\begin{remark} Relations between  dominating sets of certain graphs and graph dynamical systems previously established in \cite[Thm. 3.5, Cor. 4.7]{aledo2022fixed}, where the number of fixed points in a certain parallel sequential dynamical system and the number of a special kind of dominating sets of a bipartite graph coincide. 	\end{remark}

An immediate corollary of Theorem \ref{thm:main} (combined with Lemma \ref{prop:surjisforwardinvariant}) is the following.

\begin{corollary}
	Let $G$ and $F$ be as in Theorem \ref{thm:main}. Then for every  $s\in {\Surj}{(\bbF_2^{V})}$, there exists $n\in \bbN$ such that both $F^n(s)^{-1}(0)$ and  $F^n(s)^{-1}(1)$ are dominating sets.  
\end{corollary}
The proof follows since $F$ restricted to ${\Surj}{(\bbF_2^{V})}$ induces an $\bbN$-action and every point in a finite $\bbN$-set is eventually periodic (here fixed points assumed to have period $1$).

Let us call two configurations $s_1,s_2\in \bbF_2^V$  \emph{disjoint} if there is no such $v\in V$ with $s_1(v)=s_2(v)=1$; i.e.,  $s_1^{-1}(1)$ and $s_2^{-1}(1)$ are disjoint.
\begin{corollary}
Let $G$ and $F$ be as in Theorem \ref{thm:main}. Then there exists a domatic $n$-partition of $G$ if and only if  there exists $n$ pairwise disjoint periodic points $s_1,s_2,\dots, s_n\in {\Surj}{(\bbF_2^{V})}$  of $F$ such that $\sum_{i=1}^n s_i=c_1$.
\end{corollary}
\begin{proof}
	Suppose that $D=\{D_1,D_2,\dots, D_n\}$ is a domatic-$n$-partition of $G$. For each $i=1,\dots,n$ define $s_i:V\to \bbF_2$ by $s_i(v)=1$ if and only if $v\in D_i$. Observe that $\cup_{j\neq i}D_j$ is a dominating set of $G$, so that the pair $\{D_i,\cup_{j\neq i}D_j\}$ form a domatic $2$-partition for $G$ for every $i=1,\dots,n$. Thus, each $s_i$ is a periodic point of $F$. Since $D$ is a partition, for every $v\in V$ there exists $i\in \{1,\dots,n\}$ with $v\in D_i$, so that $s_i(v)=1$ and $s_j(v)=0$ whenever $i\neq j$. Therefore,  $(\sum_{i=1}^ns_i)(v)=1$ for all $v\in V$; i.e.,  $\sum_{i=1}^ns_i=c_1$. 
	
	Conversely, assume that  $s_1,s_2,\dots, s_n\in {\Surj}{(\bbF_2^{V})}$  are periodic points of $F$ such that $\sum_{i}s_i=c_1$. Then  $\{s_i^{-1}(1)\ \mid\ i=1,\dots,n\}$ form a domatic-$n$-partition for $G$ as each $s_i^{-1}(1)$ is a dominating set,  $s_i^{-1}(1)\cap s_j^{-1}(1)=\emptyset$ whenever $i\neq j$ and $\cup_{i}s_i^{-1}(1)=V$ since  $\sum_{i=1}^ns_i=c_1$.
\end{proof}

\begin{example}\label{ex1}
	Consider the  graph $G=\{V,E\}$ with $V=\{0,1,2,3,4\}$ and $E=\{01,12,23,20,34\}$; see the figure below
	$$\begin{tikzpicture}[node distance=2cm, font=\Large, every node/.style={scale=0.75},
		one/.style={fill=lightgray}]
		\node	 (4)      {$v_4$};
		\node        (3) [left of=4] {$v_3$};
		\node         (0) [above left of=3] {$v_0$};
		\node         (2) [below left of=3] {$v_2$};
		\node        (1) [above left of=2] {$v_1$};

		\path[-] (0)	edge             node {} (1)
		(1) edge  node {} (2)
		(2) edge  node {} (3)
		(3) edge  node {} (0)
		(4) edge  node {} (3);
	\end{tikzpicture}$$ 
	For $j\in \{0,1,2,3, 4\}$ let $i_j\in \{0,1\}$ and $s_{i_0i_1i_2i_3i_4}:V\to \bbF_2$ denote the configuration of $V$ given by $s_{i_0i_1i_2i_3i_4}(v_j)=i_j$. Then $$\bbF_2^V=\{s_{i_0i_1i_2i_3i_4} \ \mid  \  0\leq j\leq 4,\ i_j\in \bbF_2\}.$$
	The graph dynamical system $F$ associated to this graph can be viewed as follows:
	$$\begin{tikzpicture}[node distance=2cm, font=\Large, every node/.style={scale=0.85},
		one/.style={fill=lightgray},auto]

		\node[draw]          		(8) 			 			{$s_{01000}$};
		\node[draw,one]       		(21)	[right of=8] 		{$s_{10101}$};
		\node    					(23) 	[below of=8] 		{$s_{10111}$};
		\node	        		   	(7) 	[above of=8]		{$s_{00111}$};		
		\node[draw] 		  	 	(16)	[left of=8]   		{$s_{10000}$};
		\node	 	  				(15)   	[left of=16] 		{$s_{01111}$};

		\node[draw,one] 			(10) 	[right  of=21] 		{$s_{01010}$};
		\node[draw]   				(5) 	[right of=10] 		{$s_{00101}$};
		\node[draw] 		      	(2) 	[right of=5] 		{$s_{00010}$};
		\node   		     		(29) 	[right of=2] 		{$s_{11101}$};
		
		\node[draw] 	       		(1) 	[above of=5] 		{$s_{00001}$};
		\node         				(30) 	[above of=1] 				{$s_{11110}$};		
		
		\node[draw]         		(4)		[above of=10] 				{$s_{00100}$};
		\node 						(27) 	[above of=4] 			{$s_{11011}$};

		\node[draw]           		(17) 	[below of=5] 		{$s_{10001}$}; 
		\node       				(14) 	[below of=17] 				{$s_{01110}$};

		\node[draw]        			(20)	[below of=10] 				{$s_{10100}$};
		\node          				(11) 	[below of=20] 		{$s_{01011}$};

		\node		       	 		(24) 	[below of=2] 		{$s_{11000}$};
		\node          				(26) 	[above of=2] 		{$s_{11010}$};
		
		\node        				(28)	[left of=27] 	{$s_{11100}$};
		\node[draw]          		(3) 	[left of=28] 		{$s_{00011}$};
		\node[draw,one]           	(12) 	[left of=3] 	{$s_{01100}$}; 
		\node[draw,one]        		(19) 	[left of=12] 		{$s_{10011}$};
		
		\node       				(25)	[left of=11] 	{$s_{11001}$};
		\node[draw]          		(6) 	[left of=25] 		{$s_{00110}$};
		\node[draw,one]           	(9) 	[left of=6] 	{$s_{01001}$}; 
		\node[draw,one]        		(22) 	[left of=9] 		{$s_{10110}$};
		
		\node[draw,one]        		(13)	[left of=23] 	{$s_{01101}$};
		\node[draw,one]          	(18) 	[left of=13] 		{$s_{10010}$};
		
		\node[draw,one]        		(0) 	[left of=7] 		{$s_{00000}$};
		\node         				(31) 	[left of=0] 	{$s_{11111}$};

		\path[->] 
		(15) edge  node {} (16)
		(16) edge  node {} (8)
		(23) edge  node {} (8)
		(7) edge  node {} (8)
		(8) edge  node {} (21)
		
		(24) edge  node {} (5)
		(30) edge  node {} (1)
		
		(1) edge  node {} (10)
		(5) edge  node {} (10) 
		(4) edge  node {} (10) 
		
		(29) edge  node {} (2)
		(2) edge  node {} (5) 
		
		(26) edge  node {} (5)

		(27) edge  node {} (4) 
		(4) edge  node {} (10) 
		
		(28) edge  node {} (3) 
		(3) edge  node {} (12) 	
		
		(25) edge  node {} (6) 
		(6) edge  node {} (9)

		(21) edge [bend left] node {} (10)
		(10) edge [bend left] node {} (21)	
		
		(12) edge [bend left] node {} (19)
		(19) edge [bend left] node {} (12)
		
		(22) edge [bend left] node {} (9)
		(9) edge [bend left] node {} (22)

		(20) edge node {} (10)
		(11) edge node {} (20)
		
		(18) edge [bend left] node {} (13)
		(13) edge [bend left] node {} (18)	
		
		(17) edge node {} (10)
		(14) edge node {} (17)				
		
		(31) edge  node {} (0)
		(0) edge  [loop below]  node {} (0)	;
	\end{tikzpicture}
	$$	
	
	The rectangular nodes are in one-to-one correspondence with dominating sets of $G,$ where the correspondence is given by sending the element to the  preimage of $0$ under that element, see Corollary \ref{cor:inverseimageofzeroisdomset}. In this case, $G$ has $19$ dominating sets. The shaded nodes are the periodic and fixed points.  Note that the only fixed point is $s_{00000},$ the constant function at $0.$ The set of all non-zero elements inside shaded nodes admits an involution  induced by $F$ and elements in the set of orbits under this involution are in  one-to-one correspondence with  domatic $2$-partitions of $G$.  Since there is no pairwise disjoint triple in the set of periodic points, there is no domatic-$3$-partition of $G$.
\end{example}

\section{Extension to semigroup dynamical systems}\label{ssect:moredominating}
For any set $K,$ denote by $\cF(K)$ the free semigroup (with identity) generated by elements of $K.$ Let $\cX=\{F_t:S^V\to S^V \ \mid  \ t\in K\}$ be a collection of graph dynamical systems on the same dependency graph $G=(V,E).$ For $t\in K,$  one defines the action of $t\in \cF(K)$ (when considered as a word with single letter) on the configuration $s\in  S^V$ by $t\cdot s=F_t(s).$ Using this, we can induce an action of $\cF(K)$ on $S^V$ as follows. Let $w=t_{i_1}t_{i_2}\dots t_{i_n} \in \cF(K)$ and $s\in  S^V.$ Then $w\cdot s=F_{i_1}\circ F_{i_2}\dots F_{i_n}(s),$ and  $e\cdot s=s$ where $e$ is the empty word. Note that for every $w_1,w_2\in \cF(K),$ we have $w_1\cdot(w_2\cdot s)=(w_1w_2)\cdot s;$ i.e., the action is well-defined. We call such an action \emph{a $\cF(K)$-dynamical system with dependency graph $G=(V,E)$} with the state set $S.$

Suppose that we are given such an $\cF(K)$-dynamical system. Assume that both $V$ and $S$ are finite sets. For each $s\in S^V$ let $\Orb(s)=\{w\cdot s\ \mid  \ w\in\cF(K)\};$ that is, the orbit of $s$ under the action of $\cF(K).$ We abuse notation and denote the restriction of $F_t$ to $\Orb(s)$ again  by $F_t.$ Define $\ell_{\cX}(S^V)\subseteq S^V$ as follows: 
$$\ell_{\cX}(S^V)=\{s\in S^V \ \mid  \ \forall t\in K,\ F_t:\Orb(s)\to \Orb(s)\ \text{ is a permutation}\}.$$
Note that the set $\ell_{\cX}(S^V)$ can be identified with the set $A^\ell$ in \cite[Sec. 4.4, Lem. 3]{erdal2018semigroup} for the monoid $\cF(K)$ and the $\cF(K)$-set $S^V$. Note that this set is invariant under the action of $\cF(K)$; that is, for every $w\in \cF(K)$ and $s\in\ell_{\cX}(S^V),$ we have $w\cdot s\in \ell_{\cX}(S^V).$ In fact, since each $F_t$ restricted to $\ell_{\cX}(S^V)$ is a permutation, the map $$\mu_w:\Orb(s)\to \Orb(s): x\mapsto w\cdot x$$ is also a permutation, and since $S^V$ is finite and $s\in \Orb(s),$ there exists $n\in \bbZ^+$ with $w^n\cdot s=s.$  Therefore, $s\in \Orb(w\cdot s);$ i.e., $\Orb(w\cdot s)=\Orb(s).$

Observe that $\cF(K)$ acts on $\ell_{\cX}(S^V)$ by bijections.  Moreover, if $s\in \ell_{\cX}(S^V),$ then $s$ is a periodic or fixed point for every map $F_t.$ More precisely, any subset of $S^V$ that is invariant under the $\cF(K)$-action is also contained in $\ell_{\cX}(S^V)$. However, a point that is periodic or fixed in every system $F_t$ for $k\in K$ does not imply that it is in $\ell_{\cX}(S^V),$ as points in $\ell_{\cX}(S^V)$ are also periodic or fixed for compositions of such systems.

\subsubsection{Applications on independent domination}
A subset $A\subseteq V$ is called an independent set of $G=(V,E)$ if no two vertices in $A$ is adjacent. An independent set that is also dominating is called an independent dominating set. These subsets are also called  maximal independent sets as such a set can not be properly contained in a larger independent set.
Let $\cF(a,b)$ be the free semigroup with identity (i.e., free monoid) with two generators $a,b.$ 

We consider two different actions of $\cF(a,b)$ on $\bbF_2^{V}.$ Let $F_0:\bbF_2^{V}\to \bbF_2^{V}$ be the graph dynamical system given by $F_0(s)=sP_s$; so that, for every $v\in V$ $$F_0(s)(v)=s(v)\prod_{x\in N(v)} (1+ s(x)).$$ Let $F$ is as in Theorem \ref{thm:main}.  First consider the following action of $\cF(a,b)$ on $\bbF_2^{V}$ defined on generators by $a\cdot s = F^2(s)$ and let $b\cdot  s = F_0(s).$ This is the action generated by the family of systems $\cX_1=\{F^2, F_0\}.$  

It is known that any independent dominating set must be minimal, so that its complement is also a dominating set; i.e., an independent dominating set with its complement form a domatic-$2$-partition. Thus, the following corollary follows from Theorem \ref{thm:main}.
\begin{corollary}\label{cor:fixeddom}
	If  $s^{-1}(1)$ is an independent dominating set of $G,$ then $a\cdot s=s.$ 
\end{corollary}	
We have the following lemma:
\begin{lemma}\label{lem:independent}
	For $s\in \bbF_2^{V},$  we have $b\cdot s=s$ if and only if $s^{-1}(1)$ is an independent set of $G.$
\end{lemma}	
\begin{proof}
	Assume $s^{-1}(1)$ is an independent set and let $v\in V.$ If $s(v)=1,$ then  $s(x)=0$ for every $x\in N(v).$ Thus, $$(b\cdot s)(v)=s(v)\prod_{x\in N(v)}(1+ s(x))=1.$$ If $s(v)=0,$ then $(b\cdot s)(v)=0.$ Thus,   $b\cdot s=s.$ Conversely, assume $s^{-1}(1)$ is not an independent set. Then there exists $v,w\in V$ such that ${vw}\in E$ and $s(v)=s(w)=1.$ This implies $$\prod_{x\in N(v)}(1+ s(x))=0$$ as $w\in N(v)$; and thus, $(b\cdot s)(v)=0$. In other words, we have $b\cdot s\neq s.$
\end{proof}
We have the following corollary.
\begin{corollary}\label{cor:maxindep}
	For $s$ in $\bbF_2^{V}$ with $s\neq c_0,$   $s^{-1}(1)$ is a maximal independent set if and only if $s\in\ell_{\cX_1}(\bbF_2^V).$
\end{corollary}
\begin{proof}
	Assume $s^{-1}(1)$ is a maximal independent set of $G.$ Then by Lemma \ref{lem:independent} $b\cdot s=s.$ Since any maximal independent set is also a dominating set, by Corollary \ref{cor:fixeddom}, $a\cdot s=s.$ This means $s$ is a fixed point of the action; i.e.,  both $F^2$ and $F_0$ are identity maps when restricted to  $\Orb(s)=\{s\}.$ Therefore,  $s\in\ell_{\cX_1}(S^V).$

	Conversely, assume  $s\in\ell_{\cX_1}(S^V)$ with $s\neq c_0.$ Then $F^2$ and $F_0$ are permutations of  $\Orb(s).$ By Theorem \ref{thm:main} $s$ must be a fixed point of $F^2$ and both ${s}^{-1}(0)$ and ${s}^{-1}(1)$ are dominating sets of $G.$ Besides, as $\Orb(s)$ is finite, there exists $k\in \bbN^+$ with $b^k\cdot s=s.$ Note that if $s(v)=0,$ then $(b\cdot s)(v)=0$ and so $(b^k\cdot s)(v)=0.$ Hence, if $b^k\cdot s=s,$ then $k=1$; i.e.,  $b\cdot s=s.$ By Lemma \ref{lem:independent} this implies that  ${s}^{-1}(1)$ is an independent set. Therefore, we need ${s}^{-1}(1)$ to be simultaneously a dominating set and independent set. But this means ${s}^{-1}(1)$ is a independent dominating set, or equivalently, a maximal independent set.
\end{proof}

An idomatic partition of a graph $G=(V,E)$ is a domatic partition where each class in the partition is also an independent set; i.e., a partition of $V$ into independent dominating sets (see \cite[13.3.2]{zelinka2017domatic}. If number of classes in such a partition is $n,$ then we call it an idomatic-$n$-partition.

The second action of $\cF(a,b)$ on $\bbF_2^V$ that we consider is the action associated to the family of systems $\cX_2=\{F, F_0\},$ where $F$ and $F_0$ is as before. That is, for every $s\in \bbF_2^V,$ $a\cdot s=F(s)$ and $b\cdot s=F_0(s).$ Then we have the following proposition. 
\begin{corollary}\label{cor:idomatic}
	The graph $G$ has an idomatic-$2$-partition if and only if $\ell_{\cX_2}(\bbF_2^V)$ contains a configuration other than $c_0.$ If there exists $s\in  \ell_{\cX_2}(\bbF_2^V)$ with $s\neq c_0,$ then  $\{s^{-1}(1), s^{-1}(0)\}$ forms a idomatic-$2$-partition of $G,$ and any  idomatic-$2$-partision of $G$ is of this form for some $s\in\ell_{\cX_2}(\bbF_2^V).$
\end{corollary}
\begin{proof}
	Suppose that $\{D,B\}$ is an idomatic-$2$-partition of $G.$ Define $s\in \bbF_2^V$ by $s(v)=1$ if and only if $v\in D.$ By Theorem \ref{thm:main}, as $\{s^{-1}(1), s^{-1}(0)\}$ is a domatic-$2$-partition, $F(s)=1+s$ and $F(1+s)=s$. Besides, as both $s^{-1}(1)$ and $(1+s)^{-1}(1)$ are maximal independent sets, $F_0(s)=s$ and $F_0(1+s)=1+s$  by Lemma \ref{cor:maxindep}. Therefore, $s$ and $ 1+s$ are configurations in $\ell_{\cX_2}(\bbF_2^V).$ Conversely, if  $s\in  \ell_{\cX_2}(\bbF_2^V)$ and $s\neq c_0,$ then by definition of $\ell_{\cX_2}(\bbF_2^V),$ $s$ is a periodic or fixed point of both $F$ and $F_0.$ Since $s\neq c_0,$ $s$ is a periodic point of $F$ with period $2.$  Thus, $a\cdot s=F(s)=1+s$ and  $\{s^{-1}(1), (1+s)^{-1}(1)\}$ is a domatic-$2$-partition of $G$. Note that $\ell_{\cX_2}(\bbF_2^V)$ is invariant under the $\cF(a,b)$-action, so that $1+s \in \ell_{\cX_2}(\bbF_2^V)$. Besides, due to the proof of Corollary \ref{cor:maxindep}, every periodic point of $F_0$ is a fix point. Thus, both $s$ and $1+s$ are fixed points of $F_0.$ By Corollary  \ref{cor:maxindep}, both $s^{-1}(1)$ and $(1+s)^{-1}(1)=s^{-1}(0)$ are maximal independent sets and together they form an idomatic-$2$-partition of $G.$
\end{proof}
\begin{example}
	For the graph $G$ in Example \ref{ex1}, the fixed points of $F_0$ among the periodic points of $F$ are $s_{10101},\ s_{01010}$ and $c_0$. Therefore,  $$\ell_{\cX_2}(\bbF_2^V)=\{s_0,s_{10101}, s_{01010}\},$$
 so that $G$ has an idomatic-$2$-partition $\{\{v_1,v_3\},\{v_0,v_2
 v_4\}\}$. \end{example}


\begin{thebibliography}{10}
	
	\bibitem{aledo2022fixed}
	Juan~A. Aledo, Ali Barzanouni, Ghazaleh Malekbala, Leila Sharifan, and Jose~C.
	Valverde.
	\newblock Fixed points in generalized parallel and sequential dynamical systems
	induced by a minterm or maxterm {Boolean} functions.
	\newblock {\em J. Comput. Appl. Math.}, 408:13, 2022.
	\newblock Id/No 114070.
	
	\bibitem{barrett2000elements}
	C.~L. Barrett, H.~S. Mortveit, and C.~M. Reidys.
	\newblock Elements of a theory of simulation. {II}: {Sequential} dynamical
	systems.
	\newblock {\em Appl. Math. Comput.}, 107(2-3):121--136, 2000.
	
	\bibitem{barrett2011modeling}
	Chris Barrett, Harry B.~III Hunt, Madhav~V. Marathe, S.~S. Ravi, Daniel~J.
	Rosenkrantz, and Richard~E. Stearns.
	\newblock Modeling and analyzing social network dynamics using stochastic
	discrete graphical dynamical systems.
	\newblock {\em Theor. Comput. Sci.}, 412(30):3932--3946, 2011.
	
	
	
		\bibitem{crescenzi1995compendium}
Crescenzi, Pierluigi and Kann, Viggo
	\newblock A compendium of NP optimization problems.
	\newblock {URL: https://www.csc.kth.se/tcs/compendium/}, 1995.
	

	
	
	\bibitem{erdal2018semigroup}
	Mehmet~Akif Erdal and {\"O}zg{\"u}n {\"U}nl{\"u}.
	\newblock Semigroup actions on sets and the {Burnside} ring.
	\newblock {\em Appl. Categ. Struct.}, 26(1):7--28, 2018.
	
	\bibitem{eubank2004modelling}
	Stephen Eubank, Hasan Guclu, VS~Anil~Kumar, Madhav~V Marathe, Aravind
	Srinivasan, Zoltan Toroczkai, and Nan Wang.
	\newblock Modelling disease outbreaks in realistic urban social networks.
	\newblock {\em Nature}, 429(6988):180--184, 2004.
	
	\bibitem{haynes2020topic}
	Teresa~W. Haynes, Stephen~T. Hedetniemi, and Michael~A. Henning, editors.
	\newblock {\em Topics in domination in graphs}, volume~64 of {\em Dev. Math.}
	\newblock Cham: Springer, 2020.
	
	\bibitem{haynes1998fund}
	Teresa~W. Haynes, Stephen~T. Hedetniemi, and Peter~J. Slater.
	\newblock {\em Fundamentals of domination in graphs}, volume 208 of {\em Pure
		Appl. Math., Marcel Dekker}.
	\newblock New York, NY: Marcel Dekker, Inc., 1998.
	
	\bibitem{saadatpour2010attractor}
	Assieh Saadatpour, Istv{\'a}n Albert, and R{\'e}ka Albert.
	\newblock Attractor analysis of asynchronous {Boolean} models of signal
	transduction networks.
	\newblock {\em J. Theor. Biol.}, 266(4):641--656, 2010.
	
	\bibitem{sutner1989cellular}
	K.~Sutner.
	\newblock Linear Cellular Automata and
	the Garden-of-Eden.
	\newblock {\em Math. Intell.}, 11(2):49--53, 1989.
	\bibitem{zelinka2017domatic}
	Bohdan Zelinka.
	\newblock Domatic numbers of graphs and their variants: a survey.
	\newblock In {\em Domination in Graphs Advanced Topics}, pages 351--378.
	Routledge, 2017.
	
\end{thebibliography}
\end{document}